\magnification=1200

\def\A{{\cal A}}
\def\I{{\cal I}}
\def\J{{\cal J}}
\def\N{{\bf N}}
\def\Q{{\cal Q}}
\def\R{{\bf R}}
\def\U{{\cal U}}
\def\V{{\cal V}}
\def\W{{\cal W}}
\def\hal{{\vrule height 10pt width 4pt depth 0pt}}

\centerline{\bf A prime C*-algebra that is not primitive}
\medskip

\centerline{Nik Weaver\footnote{*}{Supported by NSF grant DMS-0070634\hfill\break
Math Subject Classification numbers: 46L05, 16G99}}
\bigskip
\bigskip

{\narrower{
\it We construct a non-separable C*-algebra that is prime but not primitive.}
\bigskip}
\bigskip

A C*-algebra is {\it prime} if the intersection of any two nonzero
(closed, two-sided) ideals is nonzero. It is {\it primitive} if there
exists a faithful irreducible representation.
\medskip

It is fairly easy to see that any primitive C*-algebra must be prime.
Conversely, following a related result of Kaplansky [2], Dixmier [1]
showed that every separable prime C*-algebra is primitive. The question
whether, indeed, every prime C*-algebra is primitive has remained open.
We answer this question negatively by more or less explicitly constructing
a counterexample.
\medskip

The basic idea is this. Our C*-algebra $\A$ is generated by a family of
commuting projections together with a family of partial isometries which
link selected projections. $\A$ is prime because any nonzero ideal contains
a projection and any two projections lie above a pair of projections which
are linked by a partial isometry. $\A$ is not primitive because every state
on $\A$ vanishes on the ideal generated by some projection --- and thus no
irreducible representation, indeed no cyclic representation, can be faithful.
Thus, $\A$ is prime but not primitive. The two conditions place competing
demands on the set of partial isometries: it must be sufficiently abundant
to verify the first claim and sufficiently sparse to verify the second. The
demands are met simultaneously by selecting the set of partial isometries via
a process of transfinite induction.
\medskip

I wish to thank Charles Akemann for help and encouragement.
\bigskip
\bigskip

\noindent {\bf 1. The construction}
\bigskip

We adopt the conventions that $\N = \{1, 2, 3, \ldots\}$, that $\R^n =
\{a: \{1, \ldots, n\} \to \R\}$, and that $2^{\aleph_0}$ is the smallest
ordinal of cardinality $2^{\aleph_0}$. Let $X$ be the set of all sequences
of real numbers, i.e., $X = \{x: \N \to \R\}$. Give $X$ counting measure
and let $H = l^2(X)$. The C*-algebra we construct will be naturally
represented on $H$.
\medskip

For each $a \in \R^n$ define
$$X_a = \{x \in X: x(i) = a(i)\hbox{ for }1 \leq i \leq n\}$$
and let $P_a \in B(H)$ be the orthogonal projection onto $l^2(X_a)
\subset l^2(X)$. Define
$$\U_n = \{P_a: a \in \R^n\}$$
and $\U = \bigcup \U_n$.
\medskip

Next, for any $n \in \N$ and any $a,b \in \R^n$, define a bijection
$f_{a,b}: X_b \to X_a$ by
$$f_{a,b}(x)(i) = \cases{a(i)& if $i \leq n$\cr
x(i)& if $i > n$\cr}$$
for $x \in X_b$ and $i \in \N$. Observe that any sequence in $X_b$ contains
the $n$-tuple $b$ as an initial segment and $f_{a,b}$ simply replaces this
initial segment with $a$. Since $f_{a,b}$ is a bijection, composition with
$f_{a,b}$ defines an isometry from $l^2(X_a)$ onto $l^2(X_b)$. Thus let
$V_{a,b} \in B(H)$ be the partial isometry with domain projection $P_a$
and range projection $P_b$ defined by composition with $f_{a,b}$ on $X_a$.
\medskip

The following observations will be useful in the sequel. For any $P_a, P_b
\in \U$, if $a$ is an extension of $b$ then $P_a \leq P_b$, while if neither
$a$ extends $b$ nor $b$ extends $a$ then $P_aP_b = 0$. Also, $P_a = V_{a,a}$.
Lastly, an easy induction shows that any finite product of partial isometries
of the form $V_{a,b}$ is again a partial isometry of the form $V_{a,b}$, even
if the different factors involve tuples of different lengths.
\medskip

Let $\W_1 = \U \times \U$ and let $\W_2 = \{Z \subset \U: Z$ is countable$\}$.
By standard set theory, $\U$, $\W_1$, and $\W_2$ all have cardinality
$2^{\aleph_0}$. Thus, there exist bijections $\phi_i: 2^{\aleph_0} \to \W_i$
($i = 1,2$). Fix such bijections $\phi_1$ and $\phi_2$. We now construct a
transfinite sequence of natural numbers $(n_\sigma)$ and two transfinite
sequences of tuples $(a_\sigma)$ and $(b_\sigma)$, all for $\sigma <
2^{\aleph_0}$, with the following properties:
\medskip

{\narrower{

\noindent (i) $a_\sigma, b_\sigma \in \R^{n_\sigma}$,

\noindent (ii) if $\phi_1(\sigma) = (P_a, P_b)$ then $a_\sigma$ and $b_\sigma$
properly extend $a$ and $b$, respectively, and

\noindent (iii) $a_\sigma(n_\sigma), b_\sigma(n_\sigma)
\not\in R_\sigma = S_\sigma \cup T_\sigma$ where
$$\eqalign{S_\sigma &= \{a_\tau(n_\sigma), b_\tau(n_\sigma):
\tau < \sigma\hbox{ and }n_\sigma \leq n_\tau\}\cr
T_\sigma &= \{P_c(n_\sigma): c \in \phi_2(\tau)\hbox{ for some }\tau < \sigma
\hbox{ and }n_\sigma \leq {\rm length}(c)\}\cr}$$}}

\noindent for all $\sigma < 2^{\aleph_0}$. The construction proceeds by
transfinite induction. Fix $\sigma < 2^{\aleph_0}$ and suppose $a_\tau$ 
and $b_\tau$ are defined for all $\tau < \sigma$. Say $\phi_1(\sigma)
= (P_a, P_b)$. Let $n_a$ and $n_b$ be the lengths of $a$ and $b$,
respectively, and let $n_\sigma = {\rm max}(n_a, n_b) + 1$. Then define
$R_\sigma$, $S_\sigma$, and $T_\sigma$ as in condition (iii) above. Both
$S_\sigma$ and $T_\sigma$ are subsets of $\R$ of cardinality
strictly less than $2^{\aleph_0}$, so we can find $r \in \R$ not
in $R_\sigma = S_\sigma \cup T_\sigma$. Let $a_\sigma$ and $b_\sigma$ be
$n_\sigma$-tuples which extend $a$ and $b$, respectively, and satisfy
$a_\sigma(n_\sigma) = b_\sigma(n_\sigma) = r$. This completes the definition
of $n_\sigma$, $a_\sigma$, and $b_\sigma$.
\medskip

For any $\sigma < 2^{\aleph_0}$ define $V_\sigma = V_{a_\sigma, b_\sigma}$;
then let $\V = \{V_\sigma: \sigma < 2^{\aleph_0}\}$. Let $\A$ be the
C*-algebra generated by $\U$ and $\V$.
\medskip

This completes the construction. Although a few arbitrary choices were made,
if one wanted to make the construction more concrete a definite rule could
be given for each choice, modulo the identification of $2^{\aleph_0}$ with
an ordinal. Also note that although $\A$ is nonunital, we could give it a
unit by defining $\U_0 = \{I\}$ and letting $\U = \bigcup_{n \geq 0} \U_n$.
This would not affect any part of our argument in Section 2.
\bigskip
\bigskip

\noindent {\bf 2. Proof of the result}
\bigskip

Retain the notation of the preceding section.
\bigskip

\noindent {\bf Lemma 1.} {\it The C*-algebra $\A$ is prime.}
\medskip

\noindent {\it Proof.} For any $B \in \A$ let $g_B$ be the function on $X$
defined by $g_B(x) = \langle B\chi_x, \chi_x\rangle$, where $\chi_x \in H$
is the characteristic function of $x$.
\medskip

Let $\I$ be a nonzero closed, two-sided ideal of
$\A$. Then $\I$ contains a nonzero positive operator $A$. Since $\A$ is
generated by $\U$ and $\V$, $A$ is approximated by a sequence of
$*$-polynomials in elements of $\U$ and $\V$, which means that $A$ is
in the C*-subalgebra generated by a countable subset $\Q$ of $\U \cup \V$.
\medskip

Since $A > 0$, there exists $x \in X$ such that $g_A(x) > 0$. Fix such an
$x$. Then let $B$ be a $*$-polynomial in elements of $\Q$ such that $\|B - A\|
< g_A(x)/2$. This implies that $\|g_B - g_A\|_\infty < g_A(x)/2$, which in
particular yields $|g_B(x)| > g_A(x)/2$.
\medskip

Now $B$ is a $*$-polynomial in a finite number of elements $B_i \in \Q$, so
we can find $n \in \N$ which strictly exceeds the length of any tuple $a$ such
that $B_i = P_a$ for some $i$, as well as the length of any tuples $a$ and $b$
such that $B_i = V_{a,b}$ for some $i$. Let $\alpha' = (x(1), \ldots, x(n-1))$.
Then any $*$-monomial $B'$ in the $B_i$'s is equal to some $V_{a,b}$ with
${\rm length}(a) = {\rm length}(b) \leq n-1$, and therefore $g_{B'} \equiv 0$
or $g_{B'} \equiv 1$ on $X_{\alpha'}$. Taking linear combinations shows that
$g_B$ is also constant on $X_{\alpha'}$. Moreover, since $x \in X_{\alpha'}$
and $|g_B(x)| > g_A(x)/2$, it follows that $|g_B| > g_A(x)/2$ on $X_{\alpha'}$,
and this implies that $g_A > 0$ on $X_{\alpha'}$.
\medskip

Next, since $\Q$ is countable we can find $r \in \R$ such that
$a(n) \neq r$ whenever $P_a \in \Q$ and $n \leq {\rm length}(a)$, and
$a(n), b(n) \neq r$ whenever $V_{a,b} \in \Q$ and $n \leq {\rm length}(a)
= {\rm length}(b)$. Let $\alpha$ be the $n$-tuple which extends $\alpha'$
and satisfies $\alpha(n) = r$. We claim that for any element $C$ of the
C*-algebra generated by $\Q$, the function $g_C$ is constant on $X_\alpha$ and
$\langle C\chi_y, \chi_z\rangle = 0$ for all distinct $y,z \in X_\alpha$.
If $C$ is a $*$-monomial in elements of $\Q$ then $C$ is a partial isometry,
$C = V_{a', b'}$, with $a'(n), b'(n) \neq r$ or $a'(n), b'(n)$ not defined,
either of which implies the claim. The claim then holds for $*$-polynomials
by linearity and in general by continuity.
\medskip

In particular, the claim holds for $C = A$, so that $P_\alpha AP_\alpha$ is a
scalar multiple of $P_\alpha$. But we showed above that $g_A > 0$ on
$X_{\alpha'} \supset X_\alpha$, so in fact $P_\alpha AP_\alpha$ is a
nonzero scalar multiple of $P_\alpha$, and hence $P_\alpha \in \I$.
\medskip

Let $\J$ be any other nonzero closed, two-sided ideal of $\A$. By similar
reasoning we can find some $P_\beta \in \J$. Now the pair $(P_\alpha, P_\beta)$
belongs to $\W_1$, so we have $\phi_1(\sigma) = (P_\alpha, P_\beta)$ for some
$\sigma < 2^{\aleph_0}$. Let $a_\sigma$ and $b_\sigma$ be as in Section 1.
Then $a_\sigma$ and $b_\sigma$ properly extend $\alpha$ and $\beta$,
respectively, so $P_{a_\sigma} < P_\alpha$ and $P_{b_\sigma} < P_\beta$, and
therefore $P_{a_\sigma} \in \I$ and $P_{b_\sigma} \in \J$. Furthermore, the
partial isometry $V_\sigma$ satisfies $V_\sigma P_{a_\sigma} V_\sigma^*
= P_{b_\sigma}$. Thus $P_{b_\sigma} \in \I \cap \J$. This shows that $\A$
is prime.\hfill\hal
\bigskip

\noindent {\bf Lemma 2.} {\it Let $\rho$ be a state on $\A$. Then there
is a 1-tuple $a \in \R^1$ such that $\rho$ vanishes on the ideal generated
by $P_a$.}
\medskip

\noindent {\it Proof.} For each $n \in \N$, the projections $P_b$
are orthogonal as $b$ ranges over $\R^n$. Therefore $\rho$ is nonzero
on only countably many members of each $\U_n$. It follows that the set
$$Z = \{P_b \in \U: \rho(P_b) \neq 0\}$$
is countable. Thus $Z \in \W_2$ and we have $Z = \phi_2(\sigma)$ for some
$\sigma < 2^{\aleph_0}$. Fix this value of $\sigma$, let $S$ be the set of
all values $a_\tau(1)$ and $b_\tau(1)$ such that $\tau \leq \sigma$, and
let $T$ be the set of all values $c(1)$ such that $P_c \in \phi_2(\sigma)$.
Then define a 1-tuple $a$ by letting $a(1)$ be any element of $\R$ which
avoids $S$ and $T$.
\medskip

We will show that $\rho$ vanishes on the ideal generated by $P_a$. The
main step is the following claim: for any finite product $A$ of elements
of $\U \cup \V \cup \V^*$ at least one of which is $P_a$, there exist
$n \in \N$ and $b \in \R^n$ such that
\medskip

{\narrower{
\noindent (1) $P_bA = A$,

\noindent (2) $a_\tau(n) \neq b(n)$ and $b_\tau(n) \neq b(n)$ for all
$\tau \leq \sigma$ such that $n \leq n_\tau$, and

\noindent (3) $c(n) \neq b(n)$ for all $P_c \in \phi_2(\sigma)$ such that
$n \leq {\rm length}(c)$.}
\medskip}

We will assume $A \neq 0$. In the base case, for a product of length 1 we have
$A = P_a$ and we can take $b = a$ and $n = 1$. Property (1)
is trivial, and (2) and (3) follow directly from the definition of $a$. Now
suppose the induction hypothesis holds for all products of length at most $k$
and let $A$ be a product of length $k + 1$. If the leftmost factor of $A$ is
$P_a$ then we can still take $b = a$ and $n = 1$. Otherwise, write
$A = CB$ where $C \in \U \cup \V \cup \V^*$ and $B$ is a product of length
$k$ which includes $P_a$ as a factor. Let $b$ and $n$ satisfy (1) -- (3)
for $B$. If $C \in \U$ then $b$ and $n$ also satisfy (1) -- (3) for
$A$, so assume $C \in \V \cup \V^*$. Without significant loss of generality
we can suppose $C = V_\kappa$ for some $\kappa < 2^{\aleph_0}$.
\medskip

Let $a_\kappa$ and $b_\kappa$ be as in Section 1. If $n > n_\kappa$
then ${\rm length}(a_\kappa) < {\rm length}(b)$, and
$CB = A \neq 0$ and $P_bB = B$ then imply that $b$ is an extension of
$a_\kappa$. (Otherwise the cokernel of $C = V_\kappa$, namely
$l^2(X_{a_\kappa})$, would be orthogonal to $l^2(X_b)$ and hence to
the range of $B$.) Then $C(l^2(X_b)) = l^2(X_{b'})$ where $b'$ is $b$ with
the initial $a_\kappa$ segment replaced by $b_\kappa$. Thus
${\rm length}(b') = {\rm length}(b)$ and $b'(n) = b(n)$, and $b'$
and $n$ satisfy (1) -- (3) for $A$. Property (1) holds because
$$P_{b'}A = P_{b'}CB = CP_bB = CB = A,$$
and properties (2) and (3) hold because $b'(n) = b(n)$. 
\medskip

Finally, suppose $n \leq n_\kappa$. If $\kappa \leq \sigma$ then by property
(2) we have $a_\kappa(n) \neq b(n)$. This implies that $X_{a_\kappa}$ is
disjoint from $X_b$, so that $l^2(X_{a_\kappa})$ is orthogonal to $l^2(X_b)$
and hence $A = CB = V_\kappa P_bB$ is zero just as in the last paragraph,
contradicting our assumption. So we must have $\kappa > \sigma$. Then
$b_\kappa$ and $n_\kappa$ satisfy (1) -- (3) for $A$: property (1) holds
since $P_{b_\kappa}$ is the range projection of $C$, and properties (2) and
(3) follow from property (iii) in the original definition of $a_\kappa$ and
$b_\kappa$ together with the fact that $\sigma < \kappa$. This completes the
proof of the claim.
\medskip

Now let $A$ be any finite product of elements of $\U \cup \V \cup \V^*$, at
least one of which is $P_a$, and let $b$ be as in the claim. Since
$Z = \phi_2(\sigma)$, property (3) implies that $P_b \not\in Z$, i.e.,
$\rho(P_b) = 0$. It follows from the Cauchy-Schwarz inequality for positive
functionals [3, Theorem 3.1.3] that
$$|\rho(A)|^2 = |\rho(P_bA)|^2 \leq \rho(P_b)\rho(A^*A) = 0,$$
and hence $\rho(A) = 0$. We have therefore shown that $\rho$ vanishes on any
finite product of elements of $\U \cup \U \cup \V^*$, at least one of which
is $P_a$. Taking linear combinations shows that $\rho$ vanishes on any
$*$-polynomial in elements of $\U$ and $\V$ every term of which involves
$P_a$, and by continuity $\rho$ then vanishes on the ideal generated by
$P_a$.\hfill\hal
\bigskip

\noindent {\bf Theorem.} {\it The C*-algebra $\A$ is prime but not
primitive.}
\medskip

\noindent {\it Proof.} $\A$ is prime by Lemma 1. Now let $\pi: \A \to B(K)$
be an irreducible representation on a Hilbert space $K$ and let $\xi \in K$
be a unit vector. Then $\rho(A) = \langle \pi(A)\xi, \xi\rangle$ defines a
state on $\A$. By Lemma 2, there is a nonzero ideal $\I$ on which $\rho$
vanishes. Then
$$\langle \pi(A)\pi(B)\xi, \pi(C)\xi\rangle =
\langle \pi(C^*AB)\xi, \xi\rangle = \rho(C^*AB) = 0$$
for all $A \in \I$ and $B, C \in \A$. Since $\pi$ is irreducible $\xi$ is
cyclic, so this implies that $\pi(A) = 0$ for all $A \in \I$. Thus $\pi$
is not faithful, and we conclude that $\A$ is not primitive.\hfill\hal

\bigskip
\bigskip

[1] J.\ Dixmier, Sur les C*-alg\`{e}bres, {\it Bull.\ Soc.\ Math.\
France \bf 88} (1960), 95-112.
\medskip

[2] I.\ Kaplansky, The structure of certain operator algebras, {\it Trans.\
Amer.\ Math.\ Soc.\ \bf 70} (1951), 219-255.
\medskip

[3] G.\ K.\ Pedersen, {\it C*-Algebras and Their Automorphism Groups},
Academic Press (1979).
\bigskip
\bigskip

\noindent Math Dept.

\noindent Washington University

\noindent St. Louis, MO 63130

\noindent nweaver@math.wustl.edu
\end